\documentclass[10pt]{amsart}
\usepackage{amsmath,amsfonts,amssymb,amsthm,amscd}
\newtheorem{thm}{Theorem}
\newtheorem{lem}[thm]{Lemma}
\newtheorem{defn}[thm]{Definition}
\newtheorem{prop}[thm]{Proposition}
\newtheorem{cor}[thm]{Corollary}

\newlength{\dinwidth}
\newlength{\dinmargin}
\def\F{{\mathfrak F}}

\begin{document}

\author{Liana David and Ian A.B. Strachan}

\bigskip

\title{Compatible metrics on a manifold and non-local bi-Hamiltonian structures}

\bigskip

\begin{abstract}
Given a flat metric one may generate a local Hamiltonian structure
via the fundamental result of Dubrovin and Novikov. More
generally, a flat pencil of metrics will generate a local
bi-Hamiltonian structure, and with additional quasi-homogeneity
conditions one obtains the structure of a Frobenius manifold. With
appropriate curvature conditions one may define a curved pencil of
compatible metrics and these give rise to an associated non-local
bi-Hamiltonian structure. Specific examples include the
$F$-manifolds of Hertling and Manin equipped with an invariant
metric. In this paper the geometry supporting such compatible
metrics is studied and interpreted in terms of a multiplication on
the cotangent bundle. With additional quasi-homogeneity
assumptions one arrives at a so-called weak $\F$-manifold - a
curved version of a Frobenius manifold (which is not, in general,
an $F$-manifold). A submanifold theory is also developed.
\end{abstract}

\date{\today}

\maketitle

\tableofcontents

\bigskip

\section{Introduction}

Let $M$ be a smooth manifold. The space of smooth vector fields on
$M$ will be denoted ${\mathcal X}(M)$ and the space of smooth
$1$-forms on $M$ will be denoted ${\mathcal E}^{1}(M).$ If $g$ is
a (pseudo-Riemannian) metric on $M$, we shall denote by $g^{*}$
the induced metric on $T^{*}M.$ For a vector field $X$, $g(X)$
will denote the $1$-form corresponding to $X$ and for a $1$-form
$\alpha$, $g^{*}\alpha $ will denote the vector field
corresponding to $\alpha$ (via the isomorphism defined by $g$
between $TM$ and $T^{*}M$).

In this paper we will study the geometry induced by two metrics
$g$ and $\tilde{g}$ on $M$. Unless otherwise stated, we will
always denote by $\nabla$, $R$ ($\tilde{\nabla}$, $\tilde{R}$) the
Levi-Civita connection and the curvature tensor of $g$
($\tilde{g}$ respectively). For every constant $\lambda$ let
$g_{\lambda}^{*}:=g^{*}+\lambda\tilde{g}^{*}$, which, we will
assume, will always be non-degenerate. The Levi-Civita connection
and curvature tensor of $g_{\lambda}$ will be denoted by
$\nabla^{\lambda}$ and $R^{\lambda}$ respectively.

\begin{defn}\label{def1}\cite{mok}\,.
The metrics $g$ and $\tilde{g}$ are almost compatible if the
relation
\begin{equation}\label{defalmost}
g_{\lambda}^{*}\nabla^{\lambda}_{X}\alpha = g^{*}\nabla_{X}\alpha
+\lambda\tilde{g}^{*}\tilde{\nabla}_{X}\alpha
\end{equation}
holds, for every $X\in{\mathcal X}(M)$, $\alpha\in{\mathcal E}
^{1}(M)$ and $\lambda$ constant.

\noindent The metrics $g$ and $\tilde{g}$ are compatible if they
are almost compatible and moreover the relation
$$
g_{\lambda}^{*}(R^{\lambda}_{X,Y}\alpha )= g^{*}(R_{X,Y}\alpha
)+\lambda \tilde{g}^{*}(\tilde{R}_{X,Y}\alpha )
$$
holds, for every $\alpha\in{\mathcal E}^{1}(M)$, $X,Y\in{\mathcal
X}(M)$ and $\lambda$ constant.

\noindent If $R^\lambda=0$ for all $\lambda$ then $g$, $\tilde{g}$
are said to form a flat-pencil of metrics.

\end{defn}

The motivation for this definition comes from the theory of
bi-Hamiltonian structures for equations of hydrodynamic type, i.e.
for $(1+1)$-dimensional evolution equations
\[
\frac{\partial u^i}{\partial T} = M^i_j[u^r(X,T)] \frac{\partial
u^j}{\partial X}\,.
\]
The foundational result is this area is due to Dubrovin and
Novikov \cite{DN}:

\begin{thm}\label{DNthm}
Given two functionals of hydrodynamic type (i.e. depending only on
the fields $\{u^i\}$ and not their derivatives)
\[
F=\int_{S^1} f(u)\,dX\,,\quad\quad G=\int_{S^1} g(u)\,dX
\]
the bracket
\[
\{F,G\} = \int_{S^1} \frac{\delta F}{\delta u^i} \left[ g^{ij}
\frac{d~}{dX} - g^{is} \Gamma_{sk}^j u^k_X \right] \frac{\delta
G}{\delta u^j}\,\,dX
\]
defines (in the non-degenerate case $det[g^{ij}]\neq 0$) a
Hamiltonian structure if and only if
\begin{enumerate}
\item[$\bullet$] $g^{ij}$ is symmetric, and so defines a (pseudo)-Riemannian metric;
\item[]
\item[$\bullet$] $\Gamma_{ij}^k$ are the Christoffel symbols of the Levi-Civita connection of $g\,;$
\item[]
\item[$\bullet$] the curvature tensor of $g$ is identically zero.
\end{enumerate}
Such a Hamiltonian structure is said to be of Dubrovin/Novikov
type.

\end{thm}

\noindent The concept of a bi-Hamiltonian structure was introduced
by Magri \cite{Magri}. Given two Hamiltonian structures $\{,\}_1$
and $\{,\}_2$ then one may define a new bracket
\[
\{,\}_\lambda =\{,\}_1+\lambda \{,\}_2\,.
\]
\begin{defn}
If $\{,\}_\lambda$ is a Hamiltonian structure for all $\lambda$
then the brackets $\{ ,\}_{1}$ and $\{ ,\}_{2}$ define a
bi-Hamiltonian structure.
\end{defn}
It follows immediately from the above definitions and Theorem
\ref{DNthm} that one has a bi-Hamiltonian structure of
Dubrovin/Novikov type if and only if the corresponding metrics
form a flat-pencil. The existence of such a pencil on the manifold
$M$ results in a very rich geometric structure and leads, with
various extra conditions, to $M$ being endowed with the structure
of a Frobenius manifold \cite{dubrovin}.

\begin{defn}\cite{dubrovin}\,.
$M$ is a Frobenius manifold if a structure of a Frobenius algebra (i.e. a commutative,
associative algebra with multipication denoted by $\cdot\,,$ a unity element $e$ and an
inner product $<,>$ satisfying the invariance condition $<a\cdot b,c>=<a,b\cdot c>$)
is specified on any tangent plane $T_pM$ at any point $p\in M$
smoothly depending on the point such that:
\begin{itemize}
\item[(i)] The invariant metric ${\tilde{g}}=<,>$ is a flat metric on
$M\,;$
\item[(ii)] The unity vector field $e$ is covariantly constant with respect to the Levi-Civita connection
${\tilde\nabla}$ for the metric ${\tilde{g}}$
\[
{\tilde\nabla}e=0\,;
\]
\item[(iii)] The symmetric $3$-tensor $c(X,Y,Z)={\tilde g}(X\cdot Y,Z)$ be such that the
tensor
\[
({\tilde{\nabla}}_W c)(X,Y,Z)
\]
is totally symmetric;
\item[(iv)] A vector field $E$ - the Euler vector field - must be determined on $M$ such that
\[
{\tilde\nabla}({\tilde\nabla}E)=0
\]
and that the corresponding one-parameter group of diffeomorphisms acts by conformal transformations of the metric
and by rescalings on the Frobenius algebras $T_pM\,.$
\end{itemize}
\end{defn}


Generalizations of Dubrovin/Novikov structures were introduced by
Ferapontov \cite{F}. Such structures (originally obtained by
applying the Dirac theory of constrained dynamical systems) are of
the form
\[
\{F,G\} = \int_{S^1} \frac{\delta F}{\delta u^i} \left[ g^{ij}
\frac{d~}{dX} - g^{is} \Gamma_{sk}^j u^k_X+ \sum_{\alpha}
w^i_{{\alpha}q} u^q_X \left( \nabla^\perp\right)^{-1}
w^j_{{\alpha}q}u^q_X \right] \frac{\delta G}{\delta u^j}\,\,dX
\]
where $(g,\Gamma,w)$ must satisfy certain geometric conditions,
the crucial difference being the presence of curvature. Here
the $w$ may be interpreted as Weingarteen maps and $\nabla^\perp$
as a normal connection.
Bi-Hamiltonian structures may be similarly defined. Definition
\ref{def1}, first introduced by Mokhov \cite{mok}, ensures that a
compatible pair of metrics will define a bi-Hamiltonian structure
of this generalized type, usually called a non-local
bi-hamiltonian structure. No further mention will be made in this
paper of bi-Hamiltonian structures, though this was one of the
original motivations to study the geometry of compatible metrics.

The aim of this paper is to study the geometric structures on a
manifold endowed with two compatible metrics, and conversely, to
study the geometric conditions required for two metrics to be
compatible. The constructions and results will mirror those in
\cite{dubrovin}, but it will turn out that many of the results
require only almost compatibility or compatibility and not
flatness.

The rest of the paper is laid out as follows. In Section 2 the
condition on the pair $(g,\tilde{g})$ required for almost
compatibility is derived. This condition, the vanishing of the
Nijenhuis tensor constructed from the pair, appeared in
\cite{mok}. The proof given here is shorter and coordinate free.
It is included both for completeness and to fix the notions used
in later sections. Section 3 contains the central result of the
paper: the conditions required for an almost compatible pair of
metrics to be compatible. These conditions are interpreted in
terms of an algebraic structure on the cotangent bundle in Section
4. Again, the ideas follow Dubrovin \cite{dubrovin}, but the
algebraic structure comes from compatibility, not the flatness of
the pencil. The concept of an $\F$-manifold and a weak
$\F$-manifold are introduced in Section 5. With this the
connection between (weak) quasi-homogeneous pencils of metrics and
(weak) $\F$-manifolds can be made precise. Curvature properties
are studied in Section 6. In particular, when both metrics are
flat one recovers the results of \cite{dubrovin}.

$F$-manifolds were introduced by Hertling and Manin \cite{HM}. An
application of our results is that any $F$-manifold (with Euler
field), equipped with a (non-necessarily flat) invariant metric
\[
{\tilde{g}}(a\cdot b,c)={\tilde{g}}(a,b\cdot c),
\]
will generate a pencil of compatible metrics. Thus large classes
of examples may be derived from singularity theory. However, a
weaker notion is sufficient to ensure the existence of a pencil of
compatible metrics. It was shown \cite{hert} that the $F$-manifold
condition is related to the total symmetry of the tensor
${\tilde\nabla}(\cdot).$ A weak $\F$-manifold, which ensures the
existence of compatible metrics and hence of non-local
bi-Hamiltonian structures, requires only that the tensor
${\tilde\nabla}(\cdot)(X,Y,Z,E)$ is totally symmetric in $X$, $Y$
and $Z$, where $E$ is the Euler vector field. Thus all
$F$-manifolds with compatible metrics are weak $\F$-manifolds but
not vice-versa. The different fount is used to denote the fact
that the definition of a weak $\F$-manifold includes a metric
while the definition of an $F$-manifold is metric independent. Our
results from Sections 5 and 6 can be summarized in Table 1, where the
vertical arrows denote $1:1$ correspondences.
\begin{table}
\begin{tabular}{ccccccc}
$\begin{array}{c} {\rm Frobenius} \\ {\rm manifold} \end{array}$
&$\hookrightarrow$& $\begin{array}{c} {\F-{\rm manifold~with}} \\
{\tilde{\nabla}}(e)=0\,,\mathcal{L}_E(\cdot)=\cdot\end{array}$
&$\hookrightarrow$ & $\begin{array}{c} {{\rm weak~}\F-{\rm manifold~with}} \\
{\tilde{\nabla}}(e)=0\,,\mathcal{L}_E(\cdot)=\cdot\end{array}$&$\hookrightarrow$
& weak {$\F-$manifold}\cr
 \cr
  $\updownarrow $ &&$\updownarrow $&&$\updownarrow $&& $\updownarrow $ \cr
  \cr
$\begin{array}{c} {\rm quasihomogeneous} \\ {\rm
flat~pencil}\end{array}$ &$\hookrightarrow$& $\begin{array}{c}
{\rm quasihomogeneous} \\ {\rm pencil~with}
\\{\rm curvature~condition}\end{array}$  &$\hookrightarrow$& $\begin{array}{c} {\rm quasihomogeneous}
\\ {\rm pencil} \end{array}$&$\hookrightarrow$
&$\begin{array}{c} {\rm weak}\\ {\rm quasihomogeneous} \\{\rm
pencil}\end{array}$\cr
&&&&&&\cr
\end{tabular}
\caption{}
\end{table}

The origin of this paper was one of the authors' work on the
induced structures on submanifolds of Frobenius manifolds
\cite{ian,ian2}. It is natural to consider conditions for metrics
on a submanifold, induced from a compatible pencil of metrics in
the ambient manifold, to be almost compatible and compatible. Such
questions are studied in Sections 7 and 8.

The Appendix contains a short proof that, in the semi-simple case,
almost compatibility implies compatibility, a result originally
obtained in \cite{mok}. Again, it is included here for
completeness. Such a result is of interest in the study of
semi-Hamiltonian hydrodynamic systems.

Various related results have already appeared in the literature,
but with various additional assumptions, such as semi-simplicity
or flatness of at least one of the metrics. Such distinguished
additional structures simplify many of the calculations. Here no
such simplifying assumptions are made. Finally, it should be
straightforward to extend these results to the almost Frobenius
structures introduced recently by Dubrovin \cite{dubrovin2} and
studied further by Manin \cite{manin}.

\section{Almost compatible metrics}

The following Theorem has been proved in \cite{mok}. We construct
a new shorter proof which uses a coordinate free argument.

\begin{thm}\label{ac}
The metrics $g$ and $\tilde{g}$ are almost compatible if and only
if the $(2,1)$-tensor
$$
N_{A}(X,Y):=-[A(X),A(Y)]+A[A(X),Y]+A[X,A(Y)]-A^{2}[X,Y]
$$
(with $A:TM\rightarrow TM$ defined by $A:=\tilde{g}^{*}g$)
vanishes identically.
\end{thm}

\begin{proof}
We recall that the Levi-Civita connection $\nabla$ of a metric $g$
on $M$ is determined on $TM$ by the Koszul formula: for every $X,
Y,Z\in{\mathcal X}(M)$,
\begin{align*}
g(\nabla_{Y}X,Z)&=\frac{1}{2}\{X(g(Y,Z))+Y(g(X,Z))-Z(g(X,Y))\\
&-g([X,Z],Y)-g([Y,Z],X)-g([X,Y],Z)\}.\\
\end{align*}
Let $X:=g^{*}\alpha$ and $Z:=g^{*}\gamma$. Then
$X(g(Y,Z))=g^{*}(\alpha ,d(i_{Y}\gamma )),\quad
Z(g(X,Y))=g^{*}(\gamma ,d(i_{Y}\alpha ))\,$  and
\begin{align*}
g([Y,Z],X)&=-g^{*}\left( L_{Y}(\alpha ),\gamma\right)+Y\left(g^{*}(\alpha ,\gamma )\right)\\
g([X,Y],Z)&=g^{*}\left( L_{Y}(\gamma ),\alpha\right)-
Y\left(g^{*}(\alpha ,\gamma )\right).
\end{align*}
We deduce that the metrics $g$ and $\tilde{g}$ are almost
compatible if and only if
\begin{equation}\label{almostc}
g_{\lambda}([X_{\lambda},Z_{\lambda}],Y)=g([X,Z],Y)+\lambda\tilde{g}([\tilde{X},\tilde{Z}],Y)
\end{equation}
holds, where
$$
g_{\lambda}^{*}\alpha  =X_{\lambda};\quad g^{*}\alpha =X; \quad
\tilde{g}^{*}\alpha =\tilde{X}
$$
and
$$
g_{\lambda}^{*}\gamma  =Z_{\lambda};\quad g^{*}\gamma =Z; \quad
\tilde{g}^{*}\gamma =\tilde{Z}.
$$
Since $g_{\lambda}^{*}=g^{*}+\lambda\tilde{g}^{*}$,
$X_{\lambda}=X+\lambda\tilde{X}$ and
$Z_{\lambda}=Z+\lambda\tilde{Z}.$ Relation (\ref{almostc}) becomes
equivalent with
\begin{equation}
g_{\lambda}([X+\lambda\tilde{X},Z+\lambda\tilde{Z}])=
g([X,Z])+\lambda\tilde{g}([\tilde{X},\tilde{Z}]).
\end{equation}
Note that $\tilde{X}=A(X)$ and $\tilde{Z}=A(Z).$ Applying
$g_{\lambda}^{*}$ to both terms of the above equality and
identifying the coefficients of $\lambda$ we easily obtain the
conclusion.

\end{proof}

We shall use the following Proposition in our characterization of
compatible metrics.

\begin{prop}\label{au}
Suppose the metrics $g$ and $\tilde{g}$ are almost compatible.
Then for every $\alpha,\gamma\in{\mathcal E}^{1}(M)$, the relation
\begin{equation}\label{aditional}
g^{*}(\tilde{\nabla}_{\tilde{g}^{*}\gamma}\alpha
-\nabla_{\tilde{g}^{*}\gamma}\alpha ) =
\tilde{g}^{*}(\tilde{\nabla}_{{g}^{*}\gamma}\alpha
-\nabla_{g^{*}\gamma}\alpha )
\end{equation}
holds.
\end{prop}

\begin{proof}
The Levi-Civita connection $\nabla$ on $T^{*}M$ is determined by
the formula (see also the proof of Theorem \ref{ac})
$$
2g^{*}(\nabla_{Y}\alpha ,\beta )=-g^{*}\left( i_{Y}d\beta
,\alpha\right) +g^{*}\left(i_{Y}d\alpha ,\beta\right)
+Yg^{*}(\alpha ,\beta )-g([g^{*}\alpha ,g^{*}\beta ],Y).
$$
Let $A:=\tilde{g}^{*}g$ and define $A^{*}(\alpha )(X)=\alpha (AX)$
for every $\alpha\in{\mathcal E}^{1}(M)$ and $X\in{\mathcal
X}(M)$. Then
\begin{equation}\label{usor}
\tilde{g}^{*}(\alpha ,\beta )=g^{*}(\alpha ,A^{*}\beta )
\end{equation}
and the Koszul formula for $\tilde{g}$ becomes
$$
2g^{*}(\tilde{\nabla}_{Y}\alpha ,A^{*}\beta )=-\tilde{g}^{*}\left(
i_{Y}d\beta ,\alpha\right) +\tilde{g}^{*}\left(i_{Y}d\alpha
,\beta\right) +Y\tilde{g}^{*}(\alpha ,\beta
)-\tilde{g}([\tilde{g}^{*}\alpha ,\tilde{g}^{*}\beta ],Y).
$$
It is now easy to see that
\begin{align*}
2g^{*}(\nabla_{Y}\alpha -\tilde{\nabla}_{Y}\alpha ,\beta )&=
-g^{*}(i_{Y}d\beta ,\alpha )+g^{*}(i_{Y}d(A^{-1*}\beta
),A^{*}\alpha) -g([g^{*}\alpha ,g^{*}\beta ],Y)\\
&+\tilde{g}([\tilde{g}^{*}\alpha ,\tilde{g}^{*}A^{-1*}\beta ],Y)
\end{align*}
and that
\begin{align*}
2\tilde{g}^{*}(\nabla_{Y}\alpha -\tilde{\nabla}_{Y}\alpha ,\beta
)&= -g^{*}(i_{Y}d(A^{*}\beta ),\alpha )+g^{*}(i_{Y}d\beta
,A^{*}\alpha) -g([g^{*}\alpha ,g^{*}A^{*}\beta ],Y)\\
&+\tilde{g}([\tilde{g}^{*}\alpha ,\tilde{g}^{*}\beta ],Y).
\end{align*}
In order to prove relation (\ref{aditional}) we need to show
\begin{equation}\label{ad}
g^{*}(\tilde{\nabla}_{\tilde{g}^{*}\gamma}\alpha
-\nabla_{\tilde{g}^{*}\gamma}\alpha ,\beta) =
\tilde{g}^{*}(\tilde{\nabla}_{{g}^{*}\gamma}\alpha
-\nabla_{g^{*}\gamma}\alpha ,\beta)
\end{equation}
for every $\alpha ,\beta ,\gamma\in{\mathcal E}^{1}(M)$. Let
$g^{*}\gamma =X.$ Then $\tilde{g}^{*}\gamma =A(X)$ and relation
(\ref{ad}) becomes
\begin{align*}
-(d\beta )(AX,g^{*}\alpha )&+d(A^{-1*}\beta )(AX,g^{*}A^{*}\alpha
)-g([g^{*}\alpha ,g^{*}\beta ],AX)+\tilde{g}([\tilde{g}^{*}\alpha
,\tilde{g}^{*}A^{-1*}\beta ],AX)\\
&=-d(A^{*}\beta )(X,g^{*}\alpha )+d\beta (X,g^{*}A^{*}\alpha
)-g([g^{*}\alpha ,g^{*}A^{*}\beta ],X)\\
&+\tilde{g}([\tilde{g}^{*}\alpha ,\tilde{g}^{*}\beta ],X).\\
\end{align*}
Since $\tilde{g}^{*}A^{-1*}=g^{*}$ and $g^{*}A^{*}=\tilde{g}^{*}$
(see relation (\ref{usor})), this is equivalent to
\begin{align*}
&\beta ([AX,Z]-A^{-1}[AX,AZ]-A[X,Z]+[X,AZ])\\
&-g([Z,g^{*}\beta ],AX)+\tilde{g}([AZ,g^{*}\beta ],AX)\\
&+ g([Z,\tilde{g}^{*}\beta ],X)-\tilde{g}([AZ,\tilde{g}^{*}\beta ],X)\\
&=0,
\end{align*}
where $Z:=g^{*}\alpha .$ Let $Y:=g^{*}\beta .$ The almost
compatibility property of $g$ and $\tilde{g}$ implies that the
first row of the above expression is zero and then the above
expression reduces to
\begin{equation}\label{fin}
-g([Z,Y ],AX)+\tilde{g}([AZ,Y],AX)+g([Z,AY],X)
-\tilde{g}([AZ,AY],X)=0.
\end{equation}
Using $\tilde{g}(X,AY)=\tilde{g}(Y,AX)=g(X,Y)$, relation
(\ref{fin}) becomes $\tilde{g}(N_{A}(Y,Z),X)=0$ which is obviously
true since $N_{A}=0.$

\end{proof}

\section{Compatible metrics}

\begin{thm}\label{criteriu}
Suppose the metrics $g$ and $\tilde{g}$ are almost compatible. The
following statements are equivalent:

1. The metrics $g$ and $\tilde{g}$ are compatible.

2. For every $\alpha ,\beta\in{\mathcal E}^{1}(M)$ and
$X,Y\in{\mathcal X}(M)$, the relation
\begin{equation}\label{conditiecomp2}
{g}^{*}(\tilde{\nabla}_{Y}\alpha -\nabla_{Y}\alpha ,
\tilde{\nabla}_{X}\beta -\nabla_{X}\beta )=
{g}^{*}(\tilde{\nabla}_{X}\alpha -\nabla_{X}\alpha ,
\tilde{\nabla}_{Y}\beta -\nabla_{Y}\beta )
\end{equation}
holds.

3. For every $\alpha ,\beta\in{\mathcal E}^{1}(M)$ and
$X,Y\in{\mathcal X}(M)$, the relation
\begin{equation}\label{conditiecomp}
\tilde{g}^{*}(\tilde{\nabla}_{Y}\alpha -\nabla_{Y}\alpha ,
\tilde{\nabla}_{X}\beta -\nabla_{X}\beta )=
\tilde{g}^{*}(\tilde{\nabla}_{X}\alpha -\nabla_{X}\alpha ,
\tilde{\nabla}_{Y}\beta -\nabla_{Y}\beta )
\end{equation}
holds.

\end{thm}

\begin{proof}
Note that if $h$ is a pseudo-Riemannian metric on $M$ with
Levi-Civita connection $\nabla^{\prime}$, then its curvature
$R^{\prime}$ can be written in the form
\begin{align*}
h^{*}(R^{\prime}_{h^{*}(\gamma ),X}\alpha ,\beta )&=(h^{*}\gamma
)\left( h^{*}(\nabla^{\prime}_{X}\alpha ,\beta)\right)
-X\left(h^{*}(\nabla^{\prime}_{h^{*}\gamma }\alpha
,\beta )\right)-h^{*}(\nabla^{\prime}_{[h^{*}\gamma ,X]}\alpha ,\beta )\\
&+h^{*}(\nabla^{\prime}_{h^{*}(\nabla^{\prime}_{X}\beta )}\alpha ,\gamma
)+d\alpha
(h^{*}\gamma ,h^{*}\nabla^{\prime}_{X}\beta )\\
&-h^{*}(\nabla^{\prime}_{h^{*}(\nabla^{\prime}_{X}\alpha )}\beta
,\gamma )-d\beta (h^{*}\gamma ,h^{*}\nabla^{\prime}_{X}\alpha ),
\end{align*}
where $\alpha ,\beta ,\gamma\in{\mathcal E}^{1}(M)$ and
$X\in{\mathcal X}(M).$ We shall use this observation for
$h:=g,\tilde{g},g_{\lambda}$. Identifying the coefficients of
$\lambda$ in the compatibility condition
$$
g_{\lambda}^{*}(R^{\lambda}_{g_{\lambda}^{*}(\gamma ),X}\alpha
,\beta )= g^{*}(R_{g_{\lambda}^{*}(\gamma ),X}\alpha ,\beta)+
\lambda\tilde{g}^{*}(\tilde{R}_{g_{\lambda}^{*}(\gamma ),X}\alpha,
\beta )
$$
and using relation (\ref{defalmost}), we see that $g$ and
$\tilde{g}$ are compatible if and only if the expression
\begin{align*}
E_{\alpha ,\beta ,\gamma ,X}&:=g^{*}(\nabla_{X}\alpha
,\nabla_{\tilde{g}^{*}\gamma }\beta )
+\tilde{g}^{*}(\tilde{\nabla}_{X}\alpha
,\tilde{\nabla}_{g^{*}\gamma}\beta )
-g^{*}(\nabla_{\tilde{g}^{*}\gamma}\alpha ,\nabla_{X}\beta )
-\tilde{g}^{*}(\tilde{\nabla}_{g^{*}\gamma}\alpha
,\tilde{\nabla}_{X}\beta )\\
&+\tilde{g}^{*}(\tilde{\nabla}_{X}\beta
,\nabla_{g^{*}\gamma}\alpha )
-\tilde{g}^{*}(\tilde{\nabla}_{X}\alpha
,{\nabla}_{g^{*}\gamma}\beta )
+g^{*}(\tilde{\nabla}_{\tilde{g}^{*}\gamma}\alpha ,\nabla_{X}\beta
)-g^{*}({\nabla}_{X}\alpha
,\tilde{\nabla}_{\tilde{g}^{*}\gamma}\beta )\\
\end{align*}
is zero, for every $\alpha ,\beta ,\gamma\in{\mathcal E}^{1}(M)$
and $X\in{\mathcal X}(M).$ Using Proposition \ref{au}, we notice
that
\begin{align*}
E_{\alpha ,\beta ,\gamma ,X}&= g^{*}(\nabla_{X}\alpha
,\nabla_{\tilde{g}^{*}\gamma}\beta
-\tilde{\nabla}_{\tilde{g}^{*}\gamma}\beta
)+g^{*}(\tilde{\nabla}_{\tilde{g}^{*}\gamma}\alpha-
\nabla_{\tilde{g}^{*}\gamma}\alpha ,\nabla_{X}\beta )\\
&+\tilde{g}^{*}(\tilde{\nabla}_{X}\beta
,\nabla_{g^{*}\gamma}\alpha -\tilde{\nabla}_{g^{*}\gamma}\alpha )+
\tilde{g}^{*}(\tilde{\nabla}_{X}\alpha
,\tilde{\nabla}_{g^{*}\gamma}\beta -\nabla_{g^{*}\gamma}\beta )\\
&=\tilde{g}^{*}(\tilde{\nabla}_{g^{*}\gamma}\alpha
-\nabla_{g^{*}\gamma}\alpha ,\nabla_{X}\beta
-\tilde{\nabla}_{X}\beta )\\
&+\tilde{g}^{*}(\tilde{\nabla}_{g^{*}\gamma}\beta
-\nabla_{g^{*}\gamma}\beta ,\tilde{\nabla}_{X}\alpha
-\nabla_{X}\alpha )\\
\end{align*}
and also
\begin{align*}
E_{\alpha ,\beta ,\gamma ,X}&=
g^{*}(\tilde{\nabla}_{\tilde{g}^{*}\gamma}\beta
-\nabla_{\tilde{g}^{*}\gamma}\beta ,\tilde{\nabla}_{X}\alpha
-\nabla_{X}\alpha )\\
&+g^{*}(\tilde{\nabla}_{\tilde{g}^{*}\gamma}\alpha
-\nabla_{\tilde{g}^{*}\gamma}\alpha ,{\nabla}_{X}\beta
-\tilde{\nabla}_{X}\beta ).\\
\end{align*}
The Theorem is proved.

\end{proof}

\section{Multiplication on $T^{*}M$}

In this Section we show that the (almost) compatibility condition
can also be formulated in terms of a multiplication
\lq\lq\,$\circ$" on $T^{*}M.$ This multiplication has been used
and studied in \cite{dubrovin}, when the metrics are flat. We here
extend this study to the more general case of compatible
metrics, not necessarily flat.\\

By a multiplication \lq\lq\,$\circ$" on a vector bundle $V$ we
mean a bundle map
$$
\circ :V\times V\rightarrow V.
$$
The idea of defining a multiplication  on the tangent bundle dates
back to Vaisman \cite{V}. Here a multiplication on the cotangent
bundle is required \cite{dubrovin}.

\begin{defn}\label{multipl}
An arbitrary pair of metrics $(g,\tilde{g})$ on $M$ defines a
multiplication
\begin{equation}\label{multiplication}
\alpha\circ\beta :={\nabla}_{g^{*}\alpha}(\beta
)-\tilde{\nabla}_{g^{*}\alpha}(\beta )
\end{equation}
on $T^{*}M.$
\end{defn}

Note that, in general, the multiplications determined by
$(g,\tilde{g})$ and $(\tilde{g},g)$ do not coincide. The next
Proposition is a rewriting of the relations (2.5) and (2.6) of
\cite{dubrovin}.

\begin{prop}\label{interpretare}
For every $\alpha ,\beta ,\gamma\in{\mathcal E}^{1}(M)$, the
following relation holds:
\begin{equation}\label{interpret1}
g^{*}(\alpha\circ\beta ,\gamma )=g^{*}(\alpha ,\gamma\circ\beta )
.
\end{equation}
If $g$ and $\tilde{g}$ are almost compatible, then also
\begin{equation}\label{interpret2}
\tilde{g}^{*}(\alpha\circ\beta ,\gamma )=\tilde{g}^{*}(\alpha
,\gamma\circ\beta ) .
\end{equation}
\end{prop}

\begin{proof}
Relation (\ref{interpret1}) is a consequence of the torsion free
property of the connections $\nabla$ and $\tilde{\nabla}$:
\begin{align*}
g^{*}(\alpha\circ\beta ,\gamma )-g^{*}(\alpha ,\gamma\circ\beta
)&= g^{*}({\nabla}_{g^{*}\alpha }\beta
-\tilde{\nabla}_{g^{*}\alpha }\beta ,\gamma )-
g^{*}({\nabla}_{g^{*}\gamma }\beta
-\tilde{\nabla}_{g^{*}\gamma}\beta ,\alpha )\\
&={\nabla}_{g^{*}\alpha}(\beta )\left( g^{*}\gamma\right)
-\tilde{\nabla}_{g^{*}\alpha}(\beta )\left( g^{*}\gamma\right)
-{\nabla}_{g^{*}\gamma }(\beta )\left( g^{*}\alpha\right)\\
&+\tilde{\nabla}_{g^{*}\gamma}(\beta )\left(g^{*}\alpha\right)\\
&=d\beta \left( g^{*}\alpha ,g^{*}\gamma\right) +d\beta
\left( g^{*}\gamma ,g^{*}\alpha\right)\\
&=0.\\
\end{align*}
Suppose now that $g$ and $\tilde{g}$ are almost compatible.
Relation (\ref{aditional}) of Proposition \ref{au} can be written
as
$$
\tilde{g}^{*}(\alpha\circ\beta ,\gamma
)=g^{*}(\nabla_{\tilde{g}^{*}\alpha }\beta
-\tilde{\nabla}_{\tilde{g}^{*}\alpha}\beta ,\gamma ).
$$
It follows that
\begin{align*}
\tilde{g}^{*}(\alpha\circ\beta ,\gamma )-\tilde{g}^{*}(\alpha
,\gamma\circ\beta )&={\nabla}_{\tilde{g}^{*}\alpha}(\beta ) \left(
g^{*}\gamma\right) -\tilde{\nabla}_{\tilde{g}^{*}\alpha}(\beta )
\left(
g^{*}\gamma\right) \\
&-{\nabla}_{{g}^{*}\gamma}(\beta ) \left(\tilde{
g}^{*}\alpha\right) +\tilde{\nabla}_{{g}^{*}(\gamma )}(\beta )
\left( \tilde{g}^{*}\alpha\right)\\
&=d\beta (\tilde{g}^{*}\alpha ,g^{*}\gamma )+d\beta (g^{*}\gamma
,\tilde{g}^{*}\alpha ).\\
&=0.
\end{align*}

\end{proof}

The following Proposition is a reformulation of Theorem
\ref{criteriu} and generalizes equation (2.7) of \cite{dubrovin}.

\begin{prop}
Suppose that the metrics $g$ and $\tilde{g}$ are almost
compatible. Then they are compatible if and only if the relation
\begin{equation}\label{rightsym}
(\beta\circ\gamma )\circ\alpha =(\beta\circ\alpha )\circ\gamma
\end{equation}
holds, for every $\alpha ,\beta ,\gamma\in{\mathcal E}^{1}(M)$.
\end{prop}

\begin{proof}
Relation (\ref{conditiecomp}) is equivalent with
$$
\tilde{g}^{*}(\lambda\circ\alpha, \beta\circ\gamma
)=\tilde{g}^{*}(\lambda\circ\gamma ,\beta\circ\alpha ),
$$
or, using (\ref{interpret2}), with
$$
\tilde{g}^{*}(\lambda ,(\beta\circ\gamma )\circ\alpha )
=\tilde{g}^{*}(\lambda ,(\beta\circ\alpha )\circ\gamma ).
$$
Since $\tilde{g}$ is non-degenerate, the conclusion follows.
\end{proof}

The following Lemma relates, in a nice way, the curvatures of $g$
and $\tilde{g}$ with \lq\lq\,$\circ$"-multiplication. We state it
for completeness and we leave its proof to the reader.

\begin{lem}
Let $(g,\tilde{g})$ be an arbitrary pair of metrics on $M$, with corresponding
multiplication $\circ$ on $T^{*}M.$ Then
\begin{align*}
R_{g^{*}\alpha ,g^{*}\beta}(\gamma )&=\tilde{R}_{g^{*}\alpha
,g^{*}\beta}(\gamma )+\tilde{\nabla}_{g^{*}\alpha}(\circ )(\beta
,\delta )-\tilde{\nabla}_{g^{*}\beta}(\circ )(\alpha ,\delta )\\
&+\alpha\circ (\beta\circ\delta )-(\alpha\circ\beta)\circ\delta
-\beta\circ (\alpha\circ\delta )+(\beta\circ\alpha )\circ\delta
\end{align*}
for every $\alpha ,\beta ,\gamma, \delta\in{\mathcal E}^{1}(M).$
\end{lem}

\noindent Note that if $R={\tilde R}=0$ then one may integrate the
equation
\[
\tilde{\nabla}_{g^{*}\alpha}(\circ )(\beta ,\delta
)-\tilde{\nabla}_{g^{*}\beta}(\circ )(\alpha ,\delta )=0
\]
in terms of potential functions. The remaining part of the
equation - the defining condition for a Vinberg or pre-Lie algebra
- gives a differential equation for the potential functions. The
integrability of this differential equation was established by
Ferapontov \cite{ferapontov2} and Mokhov \cite{mok}.

\section{Quasi-homogeneous pencil of metrics and weak
$\F$-manifolds}

We now turn our attention to the case of (weak) quasi-homogeneous
pencils and the parallel notion of (weak) $\F$-manifolds. The aim
of this Section is to prove the last two
vertical 1-1 correspondences in Table 1 of the
introduction.

\begin{defn}(see \cite{dubrovin})
A pair $(g,\tilde{g})$ of compatible metrics on $M$ is called a
(regular) quasi-homogeneous pencil of degree $d$ if the following
two conditions are satisfied:

\begin{enumerate}

\item There is a smooth function $f$ on $M$ such that the vector fields
$E:=\mathrm{grad}_{g}(f)$ and $e:=\mathrm{grad}_{\tilde{g}}(f)$
have the following properties:

$$
[e,E]=e;\quad L_{E}(g^{*})=(d-1)g^{*};\quad
L_{e}(g^{*})=\tilde{g}^{*};\quad L_{e}(\tilde{g}^{*})=0.
$$

\item The operator $T(u):=\frac{d-1}{2}u +u (\tilde{\nabla}E)$ is an automorphism
of $T^{*}M.$
\end{enumerate}
\end{defn}

\noindent\textbf{Remark:} The following facts hold:
\begin{enumerate}

\item Since $L_{e}(\tilde{g})=0$ and $e=\mathrm{grad}_{\tilde{g}}(f)$,
$e$ is $\tilde{\nabla}$-parallel.

\item The conditions $L_{E}({g})=(1-d){g}$ and
$E=\mathrm{grad}_{g}(f)$, easily imply that
\begin{equation}\label{cov}
\nabla_{X}(E)=\frac{1-d}{2}X,
\end{equation}
for every $X\in{\mathcal X}(M)$. Also, $E$ is a conformal Killing
vector field of the metric $\tilde{g}$:
\begin{align*}
L_{E}(\tilde{g}^{*})&=L_{E}L_{e}(g^{*})=L_{e}L_{E}(g^{*})+L_{[E,e]}(g^{*})\\
&=(d-1)L_{e}(g^{*})-L_{e}(g^{*})=(d-2)\tilde{g}^{*}.
\end{align*}

\end{enumerate}

\vspace{10pt} A consequence of relation (\ref{cov}) is the
following Proposition, which justifies Definition \ref{wp}.

\begin{prop}\label{ajutatoare}
For every $u\in{\mathcal E}^{1}(M)$, $T(u )=df\circ u .$
\end{prop}

\begin{proof}
This is just a simple calculation: for every $X\in{\mathcal
X}(M)$,
\begin{align*}
(df\circ u )(X)&=\nabla_{g^{*}(df)}(u
)(X)-\tilde{\nabla}_{g^{*}(df)}(u )(X)\\
&=\nabla_{E}(u )(X)-\tilde{\nabla}_{E}(u )(X)\\
&=\nabla_{E}(u )(X)-E(u (X)) +u
(\tilde{\nabla}_{E}X)\\
&=-u (\nabla_{X}E) +u
(\tilde{\nabla}_{X}E)\\
&=\frac{d-1}{2}u (X)+u(\tilde{\nabla}_{X}E)
\end{align*}
where in the last equality we have used relation (\ref{cov}).
\end{proof}

\begin{defn}\label{wp}
A pair of compatible metrics $(g,\tilde{g})$ is a (regular) weak
quasi-homogeneous pencil of bi-degree $(d,D)$ if the following two
conditions are satisfied:

\begin{enumerate}

\item There is a vector field $E$ on $M$ with the properties:
$L_{E}(g)=(1-d)g$, $L_{E}(\tilde{g})=D\tilde{g}$.

\item The operator
$T(u):=\frac{d-1}{2}u+u(\tilde{\nabla}E)$ is an automorphism of
$T^{*}M$ and $T(u)=g(E)\circ u$, for every $u\in{\mathcal
E}^{1}(M).$
\end{enumerate}

\end{defn}

Note that, any quasi-homogeneous pencil of degree $d$ is also weak
quasi-homogeneous of bi-degree $(d,2-d).$

\vspace{10pt}
We now introduce the parallel notions of (weak)
$\F$-manifolds:

\begin{defn}\label{F}
A manifold $M$ with a multiplication $\lq\lq~\cdot$" on its
tangent bundle, a vector field $E$ and a metric $\tilde{g}$ is
called an $\F$-manifold if the following conditions are satisfied:

\begin{enumerate}

\item  the multiplication $\lq\lq~\cdot$" is associative, commutative and has
unity $e$;

\item the vector field $E$ (called the Euler vector field) admits
an inverse $E^{-1}$ with respect to the multiplication
$\lq\lq~\cdot$", satisfies $L_{E}(\cdot )=k\cdot$,
$L_{E}(\tilde{g} )=D\tilde{g}$ and the operator
$T(u):=\frac{D+k}{2}u-\tilde{g}(\tilde{\nabla}_{\tilde{g}^{*}(u)}E)$
is an automorphism of $T^{*}M.$

\item the metric $\tilde{g}$ is $\lq\lq~\cdot$"-invariant:
$\tilde{g}(X\cdot Y,Z)=\tilde{g}(X,Y\cdot Z)$, for every
$X,Y,Z\in{\mathcal X}(M)$.

\item the $(4,0)$-tensor $\tilde{\nabla}(\cdot )$ of $M$ defined by
$$
\tilde{\nabla}(\cdot
)(X,Y,Z,V):=\tilde{g}\left(\tilde{\nabla}_{X}(\cdot
)(Y,Z),V\right)
$$
is symmetric in all its arguments.

\end{enumerate}

\end{defn}

\noindent\textbf{Remark:} By a result of Hertling \cite{hert}, all
$\F$-manifolds (originally defined in \cite{HM}) are
$F$-manifolds, i.e. the multiplication \lq\lq\,$\cdot$" satisfies
$$
L_{X\cdot Y}(\cdot )=Y\cdot L_{X}(\cdot )+X\cdot L_{Y}(\cdot )
$$
for every $X,Y\in{\mathcal X}(M)$ and also the $1$-form $g(e)$ is
closed. The different typeface is used to denote the additional
structures not present in the definition of an $F$-manifold.

\begin{defn}\label{wF}
A weak $\F$-manifold satisfies all the conditions of an
$\F$-manifold except $(4)\,,$ which is replaced by the weaker
condition:

\begin{equation}\label{s}
\tilde{\nabla}(\cdot )(X,Y,Z,E)=\tilde{\nabla}(\cdot)(E,X,Y,Z)
\end{equation}
for every $X,Y,Z\in{\mathcal X}(M),$ where $E$ is the Euler vector
field.

\end{defn}

\noindent The tensor ${\tilde{\nabla}}(\cdot)$ is automatically
symmetric in the last three variables, this following from the
invariance property of the metric.

\subsection{From weak $\F$-manifolds to (weak) quasi-homogeneous
pencils}

\vspace{10pt}

We shall now prove the following theorem.
\begin{thm}\label{main1}
Let $(M,\cdot ,\tilde{g},E)$ be a weak $\F$-manifold with
$L_{E}(\tilde{g})=D\tilde{g}$, $L_{E}(\cdot )=k\cdot$ and identity
$e$. Define the metric $g$ by $g^{*}\tilde{g}=E\cdot .$ The
following facts hold:


1. The pair $(g,\tilde{g})$ is a weak quasi-homogeneous pencil of
bi-degree $(1+k-D,D).$

2. If $e$ is $\tilde{\nabla}$-parallel and
$k=1$, then in a neighborhood of any point of $M$ the pair
$(g,\tilde{g})$ is a quasi-homogeneous pencil of degree $2-D.$ If
moreover $M$ is simply connected the pair $(g,\tilde{g})$ is a
global quasi-homogeneous pencil.

\end{thm}

We divide the proof into several steps.

\begin{lem}\label{Fac}
Let $(M,\cdot ,\tilde{g},E)$ be a weak $\F$-manifold. Then
$N_{E\cdot}=0.$
\end{lem}

\begin{proof}
The torsion free property of $\tilde{\nabla}$ implies that
\begin{align*}
N_{E\cdot}(X,Y)&=-\tilde{\nabla}_{E\cdot X}(E)\cdot
Y-\tilde{\nabla}_{E\cdot X}(\cdot )(E,Y)+\tilde{\nabla}_{E\cdot
Y}(E)\cdot X\\
&+\tilde{\nabla}_{E\cdot Y}(\cdot )(E,X)-E\cdot
X\cdot\tilde{\nabla}_{Y}(E)-E\cdot\tilde{\nabla}_{Y}(\cdot )(E,X)\\
&+E\cdot
Y\cdot\tilde{\nabla}_{X}(E)+E\cdot\tilde{\nabla}_{X}(\cdot )(E,Y),
\end{align*}
for every $X,Y\in{\mathcal X}(M).$ From relation (\ref{s}) and the
commutativity and associativity of $\lq\lq~\cdot$", we see that
\begin{align*}
N_{E\cdot}(X,Y)&=-\tilde{\nabla}_{E\cdot X}(E)\cdot
Y-\tilde{\nabla}_{E}(\cdot )(E\cdot X,Y)+\tilde{\nabla}_{E\cdot Y}(E)\cdot X\\
&+\tilde{\nabla}_{E}(\cdot )(X,E\cdot Y)-E\cdot
X\cdot\tilde{\nabla}_{Y}(E)+E\cdot Y\cdot\tilde{\nabla}_{X}(E)\\
&=L_{E}(E\cdot X)\cdot Y+E\cdot X\cdot L_{E}(Y)-X\cdot
L_{E}(E\cdot Y)-E\cdot Y\cdot L_{E}(X)\\
&=0,
\end{align*}
where in the last equality we have used $L_{E}(\cdot )=k\cdot .$

\end{proof}

\begin{prop}\label{Fman}
Let $\lq\lq~\cdot$" be an associative, commutative, with identity
multiplication on $TM$, $\tilde{g}$ a $\lq\lq~\cdot$"-invariant
metric on $M$ and $E$ a vector field on $M$ which satisfies
$L_{E}(\cdot )=k\cdot$ and $L_{E}(\tilde{g})=D\tilde{g}$, for
$D\,,k$ constant. Suppose that $E\cdot $ is an automorphism of
$TM$. Define a new metric $g$ on $M$ by $g^{*}\tilde{g}=E\cdot .$
If $E^{\flat}:=\tilde{g}(E)$, then
\begin{align*}\label{ec}
2g^{*}({\nabla}_{\tilde{g}^{*}\gamma}\alpha
-\tilde{\nabla}_{\tilde{g}^{*}\gamma}\alpha ,\beta )&=
\tilde{g}\left( (D+k)\alpha -2\tilde{\nabla}_{\tilde{g}^{*}\alpha
}(E^{\flat}),\gamma\cdot\beta\right) -\tilde{\nabla}(\cdot
)(E^{\flat},\alpha ,\gamma ,\beta )\\
&+\tilde{\nabla} (\cdot )(\gamma ,\beta ,E^{\flat},\alpha
)-\tilde{\nabla} (\cdot
)(\alpha ,\beta ,E^{\flat},\gamma )\\
&+\tilde{\nabla}(\cdot )(\beta ,E^{\flat},\alpha ,\gamma )
+\gamma\left( N_{E\cdot}(\tilde{g}^{*}\alpha ,\tilde{g}^{*}\beta
)\cdot E^{-1}\right)
\end{align*}
for every $\alpha ,\beta ,\gamma\in{\mathcal E}^{1}(M)$. In
particular,
\begin{equation}\label{ec}
{\nabla}_{\tilde{g}^{*}\gamma}\alpha
-\tilde{\nabla}_{\tilde{g}^{*}\gamma}\alpha = \frac{1}{2}\left(
(D+k)\alpha -2\tilde{\nabla}_{\tilde{g}^{*}(\alpha
)}(E^{\flat})\right)\cdot (E^{\flat})^{-1}\cdot\gamma
\end{equation}
if and only if relation (\ref{s}) is satisfied.
\end{prop}

\begin{proof}
The multiplication $\lq\lq~\cdot$" on $TM$ together with the
metric $\tilde{g}$ induce a multiplication on $T^{*}M.$ For
$Y\in{\mathcal X}(M)$, let $Y^{\flat}:=\tilde{g}(Y).$ From the
proof of Proposition \ref{au} (with $g$ and $\tilde{g}$
interchanged), we obtain
\begin{align*}
2{g}^{*}({\nabla}_{Y}\alpha -\tilde{\nabla}_{Y}\alpha ,\beta )&=
d(\beta\cdot E^{\flat})(Y,\tilde{g}^{*}\alpha )-d\beta
(Y,\tilde{g}^{*}(E^{\flat}\cdot\alpha ))\\
&+\tilde{g}(E\cdot [\tilde{g}^{*}\alpha ,\tilde{g}^{*}\beta ]
-[E\cdot\tilde{g}^{*}\alpha ,\tilde{g}^{*}\beta ],Y)\\
&+\tilde{g}\left( N_{E\cdot}(\tilde{g}^{*}\alpha
,\tilde{g}^{*}\beta ),E^{-1}\cdot Y\right) ,
\end{align*}
since in our case $A^{*}(\alpha )=\alpha\cdot E^{\flat}$. Let
$$
E_{1}:=d(\beta\cdot E^{\flat})(Y,\tilde{g}^{*}\alpha )-d\beta
\left( Y,\tilde{g}^{*}(E^{\flat}\cdot\alpha )\right)
$$
and
$$
E_{2}=\tilde{g}\left( E\cdot [\tilde{g}^{*}\alpha
,\tilde{g}^{*}\beta ] -[E\cdot\tilde{g}^{*}\alpha
,\tilde{g}^{*}\beta ],Y\right) .
$$
Then
\begin{align*}
E_{1}&=\tilde{\nabla}_{Y}(\beta\cdot
E^{\flat})(\tilde{g}^{*}\alpha
)-\tilde{\nabla}_{\tilde{g}^{*}\alpha}(\beta\cdot E^{\flat})(Y)
-\tilde{\nabla}_{Y}(\beta )\left(\tilde{g}^{*}(E^{\flat}\cdot
\alpha )\right)
+\tilde{\nabla}_{\tilde{g}^{*}(E^{\flat}\cdot\alpha )}(\beta
)(Y)\\
&=\tilde{g}^{*}\left(\beta
,\tilde{\nabla}_{Y}(E^{\flat})\cdot\alpha -
\tilde{\nabla}_{\tilde{g}^{*}\alpha}(E^{\flat})\cdot
Y^{\flat}\right) +\tilde{\nabla}_{\tilde{g}^{*} (\alpha )\cdot
E}(\beta )(Y)- \tilde{\nabla}_{\tilde{g}^{*}\alpha}(\beta )(Y\cdot
E)\\
&+\tilde{\nabla}(\cdot )(Y^{\flat},\beta ,E^{\flat},\alpha )-
\tilde{\nabla}(\cdot )(\alpha ,\beta ,E^{\flat},Y^{\flat}).\\
\end{align*}

On the other hand, since $\tilde{\nabla}$ is torsion free and
$\tilde{g}$ is \lq\lq\,$\cdot$"-invariant,
\begin{align*}
E_{2}&=\tilde{g}\left( [\tilde{g}^{*}\alpha ,\tilde{g}^{*}\beta
],E\cdot
Y\right) -\tilde{g}\left( [E\cdot\tilde{g}^{*}\alpha ,\tilde{g}^{*}\beta ],Y\right)\\
&=\tilde{\nabla}_{\tilde{g}^{*}\alpha}({\beta})(E\cdot
Y)-\tilde{\nabla}_{\tilde{g}^{*}\beta}(\alpha )(E\cdot Y)\\
&-\tilde{\nabla}_{E\cdot\tilde{g}^{*}\alpha}(\beta
)(Y)+\tilde{g}(\tilde{\nabla}_{\tilde{g}^{*}\beta}(E\cdot\tilde{g}^{*}\alpha ),Y)\\
&=\tilde{\nabla}_{\tilde{g}^{*}\alpha }(\beta )(E\cdot
Y)-\tilde{\nabla}_{E\cdot\tilde{g}^{*}\alpha}(\beta )(Y)\\
&+\left(\tilde{\nabla}_{\tilde{g}^{*}\beta}(E^{\flat})\cdot\alpha
\right) (Y)+\tilde{\nabla}(\cdot )(\beta
,E^{\flat},\alpha ,Y^{\flat})\\
&=\tilde{\nabla}_{\tilde{g}^{*}\alpha }(\beta )(E\cdot
Y)-\tilde{\nabla}_{E\cdot\tilde{g}^{*}\alpha}(\beta )(Y)\\
&+\tilde{\nabla}_{\tilde{g}^{*}\beta}(E^{\flat})\left( Y\cdot
\tilde{g}^{*}\alpha\right)+\tilde{\nabla} (\cdot )(\beta
,E^{\flat},\alpha
,Y^{\flat})\\
&=\tilde{\nabla}_{\tilde{g}^{*}\alpha}(\beta )(E\cdot
Y)-\tilde{\nabla}_{E\cdot\tilde{g}^{*}\alpha}(\beta )(Y)\\
&+d(E^{\flat})\left(\tilde{g}^{*}\beta ,Y\cdot \tilde{g}^{*}\alpha
\right)+\tilde{\nabla}_{\tilde{g}^{*}(\alpha )\cdot
Y}(E^{\flat})\left(\tilde{g}^{*}\beta \right)\\
&+\tilde{\nabla}(\cdot )(\beta ,E^{\flat},\alpha ,Y^{\flat}).
\end{align*}
We deduce that
\begin{align*}
E_{1}+E_{2}&=\tilde{g}^{*}\left(\beta
,\tilde{\nabla}_{Y}(E^{\flat})\cdot\alpha
-\tilde{\nabla}_{\tilde{g}^{*}\alpha}(E^{\flat})\cdot Y^{\flat}
-i_{\tilde{g}^{*}(\alpha )\cdot Y}d(E^{\flat})
+\tilde{\nabla}_{\tilde{g}^{*}(\alpha )\cdot
Y}(E^{\flat})\right)\\
&+\tilde{\nabla}(\cdot )(Y^{\flat},\beta ,E^{\flat},\alpha
)-\tilde{\nabla}(\cdot )(\alpha ,\beta ,E^{\flat},Y^{\flat})
+\tilde{\nabla}(\cdot )(\beta ,E^{\flat},\alpha ,Y^{\flat}).
\end{align*}

Let $\tilde{g}^{*}\alpha =X$ and $\tilde{g}^{*}\beta =Z.$ Using
$L_{E}(\tilde{g})=D\tilde{g}$, $E_{1}+E_{2}$ becomes
\begin{align*}
\tilde{g}(\tilde{\nabla}_{Y}(E)\cdot
X,Z)&-\tilde{g}(\tilde{\nabla}_{X}E,Y\cdot
Z)-\tilde{g}(\tilde{\nabla}_{X\cdot Y}E,Z)
+{g}(\tilde{\nabla}_{Z}E,X\cdot
Y)+\tilde{g}(\tilde{\nabla}_{X\cdot Y}E,Z)\\
&+\tilde{\nabla}(\cdot )(Y,Z,E,X)-\tilde{\nabla}(\cdot
)(X,Z,E,Y)+\tilde{\nabla}(\cdot )(Z,E,X,Y)\\
=D\tilde{g}(X\cdot Y, Z)&+\tilde{g}(\tilde{\nabla}_{Y}(E)\cdot
X,Z)-\tilde{g}(\tilde{\nabla}_{X}(E)\cdot
Y,Z)-\tilde{g}(\tilde{\nabla}_{X\cdot Y}E,Z)\\
&+\tilde{\nabla}(\cdot )(Z,E,X,Y)+\tilde{\nabla}(\cdot
)(Y,Z,E,X)-\tilde{\nabla}(\cdot )(X,Z,E,Y).\\
\end{align*}

If follows that
\begin{align*}
2{g}^{*}\left({\nabla}_{Y}X^{\flat}
-\tilde{\nabla}_{Y}X^{\flat},Z^{\flat}\right)&=\tilde{g}\left(
(DY+\tilde{\nabla}_{Y}E)\cdot X-\tilde{\nabla}_{X}(E)\cdot Y
-\tilde{\nabla}_{X\cdot Y}E,Z\right)\\
&+\tilde{\nabla} (\cdot )(Y,Z,E,X)-\tilde{\nabla} (\cdot
)(X,Z,E,Y)+\tilde{\nabla} (\cdot
)(Z,E,X,Y)\\
&+\tilde{g}\left( N_{E\cdot}(X,Z),E^{-1}\cdot Y\right)\\
&=\tilde{g}\left( (DY+\tilde{\nabla}_{Y}E)\cdot
X+\tilde{\nabla}_{E}(X\cdot
Y)-\tilde{\nabla}_{E}(X)\cdot Y,Z\right)\\
&+\tilde{g}\left(-\tilde{\nabla}_{E}(Y)\cdot
X-\tilde{\nabla}_{X}(E)\cdot Y-
\tilde{\nabla}_{X\cdot Y}(E),Z\right)\\
&+\tilde{\nabla} (\cdot )(Y,Z,E,X)-\tilde{\nabla} (\cdot
)(X,Z,E,Y)+\tilde{\nabla} (\cdot )(Z,E,X,Y)\\
&-\tilde{\nabla}(\cdot )(E,X,Y,Z) +
\tilde{g}\left( N_{E\cdot}(X,Z),E^{-1}\cdot Y\right) ,\\
\end{align*}
where in the second equality we have added the (identically zero)
expression
$$
\tilde{\nabla}_{E}(X\cdot Y)-\tilde{\nabla}_{E}(X)\cdot
Y-X\cdot\tilde{\nabla}_{E}(Y)-\tilde{\nabla}_{E}(\cdot )(X,Y).
$$
Using now $L_{E}(\cdot )=k\cdot$ we obtain
\begin{align*}
2{g}^{*}\left({\nabla}_{Y}X^{\flat}
-\tilde{\nabla}_{Y}X^{\flat},Z^{\flat}\right)&=\tilde{g}\left(
(DY+\tilde{\nabla}_{Y}E)\cdot X+L_{E}(X)\cdot
Y+X\cdot L_{E}(Y)+kX\cdot Y,Z\right)\\
&-\tilde{g}\left(\tilde{\nabla}_{E}(X)\cdot
Y+\tilde{\nabla}_{E}(Y)\cdot
X+\tilde{\nabla}_{X}(E)\cdot Y,Z\right)\\
&+\tilde{\nabla} (\cdot )(Y,Z,E,X)-\tilde{\nabla} (\cdot
)(X,Z,E,Y)+\tilde{\nabla} (\cdot )(Z,E,X,Y)\\
&-\tilde{\nabla}(\cdot )(E,X,Y,Z) +\tilde{g}\left(
N_{E\cdot}(X,Z),E^{-1}\cdot Y\right)\\
&=\tilde{g}\left( (D+k)X\cdot Y-2\tilde{\nabla}_{X}(E)\cdot
Y,Z\right) +\tilde{\nabla}(\cdot )(Y,Z,E,X)\\
&-\tilde{\nabla}(\cdot )(X,Z,E,Y)+\tilde{\nabla}(\cdot
)(Z,E,X,Y)-\tilde{\nabla}(\cdot )(E,X,Y,Z)\\
&+\tilde{g}\left( N_{E\cdot}(X,Z),E^{-1}\cdot Y\right) .\\
\end{align*}
The first statement of the Proposition is proved. We easily notice
that
\begin{align*}
-\tilde{\nabla}(\cdot )(E,X,Y,Z)&+\tilde{\nabla}(\cdot )
(Y,Z,E,X)-\tilde{\nabla}(\cdot )(X,Z,E,Y)+\tilde{\nabla}(\cdot
)(Z,E,X,Y)\\
+&\tilde{g}\left( N_{E\cdot}(X,Z),E^{-1}\cdot Y\right) =0
\end{align*}
for every $X,Y,Z\in{\mathcal X}(M)$ is equivalent with relation
(\ref{s}) (exchange $X$ and $Z$ in the above equality, use the
symmetry of $\tilde{\nabla}(\cdot )$ in the last three arguments
and Lemma \ref{Fac}).
\end{proof}

\textit{Proof of the theorem:} Lemma \ref{Fac} implies that the
metrics $g$ and $\tilde{g}$ are almost compatible. The
compatibility condition (\ref{conditiecomp}) is trivially
satisfied using relation (\ref{ec}) and the \lq\lq\,$\cdot$"-invariance
of $\tilde{g}.$ Since $g^{*}\tilde{g}=E\cdot$ and $E$ is the Euler
vector field,
$L_{E}(g^{*})\tilde{g}+Dg^{*}\tilde{g}=kg^{*}\tilde{g}$ or
$L_{E}(g^{*})=(k-D)g^{*}.$ Relation (\ref{ec}) can be written in
the form $g\tilde{g}^{*}(\gamma )\circ\alpha =T(\alpha )\cdot
(E^{\flat})^{-1}\cdot\gamma$, where $T$ is precisely the operator
associated to the pair $(g,\tilde{g})$ as in Definition \ref{wp}.
For $\gamma :=\tilde{g}(E)$ we obtain $g(E)\circ \alpha =T(\alpha
).$ It follows that $(g,\tilde{g})$ is a weak quasi-homogeneous
pencil of bi-degree $(1+k-D,D).$ The first statement is proved.

Suppose now that $k=1$ and $\tilde{\nabla}(e)=0.$ In particular,
$L_{e}(\tilde{g})=0$. Since $\tilde{\nabla}$ is torsion free,
$d\tilde{g}(e)=0$ and at least locally there is a smooth function
$f$ such that $\tilde{g}(e)=df.$ Since $g^{*}\tilde{g}=E\cdot$,
$E=g^{*}\tilde{g}(e)=\mathrm{grad}_{g}(f).$ Also, $[e,E]=e$
because $k=1$. In order to prove that $(g,\tilde{g})$ is a
quasi-homogeneous pencil (of degree $2-D$), we still need to show
that $L_{e}(g^{*})=\tilde{g}^{*}.$ For this, consider the equality
(which follows from $g^{*}\tilde{g}=E\cdot$)
$$
g^{*}(\alpha ,\beta )=\tilde{g}^{*}(\alpha ,\beta\cdot E^{\flat})
$$
and take its Lie derivative in the direction of $e$. From
$[e,E]=e$ and $L_{e}(\tilde{g})=0$ we get
$$
L_{e}(g^{*})(\alpha ,\beta )=\tilde{g}^{*}\left(\alpha
,L_{e}(\cdot )(\beta ,E^{\flat})\right) +\tilde{g}^{*}(\alpha
,\beta ).
$$
On the other hand, from $\tilde{\nabla}(e)=0$ and relation
(\ref{s}), we easily see that
\begin{align*}
L_{e}(\cdot )(X,E)&=-[E\cdot X,e]+[X,e]\cdot E-X\\
&=\tilde{\nabla}_{e}(E\cdot X)+[X,e]\cdot E-X\\
&=\tilde{\nabla}_{e}(E)\cdot X+\tilde{\nabla}_{E}(\cdot
)(e,X)+E\cdot\tilde{\nabla}_{e}(X)\\
&+[X,e]\cdot E-X=0,\\
\end{align*}
which implies that $L_{e}(\cdot )(\beta
,E^{\flat})=\tilde{g}\left( L_{e}(\cdot )(\tilde{g}^{*}\beta
,E)\right) =0$, because $L_{e}(\tilde{g})=0$. It follows that
$L_{e}(g^{*})=\tilde{g}^{*}$ and the Theorem is proved.

\subsection{From quasi-homogeneous pencil of metrics to weak
$\F$-manifolds}

\begin{thm}\label{cvasihom}
Let $(g,\tilde{g})$ be a weak quasi-homogeneous pencil as in
Definition \ref{wp}. Define a new multiplication $u\cdot v:=u\circ
T^{-1}(v)$ on $T^{*}M$ and denote also by $\lq\lq~\cdot$" the
induced multiplication on $TM$, using the metric $\tilde{g}.$ The
following statements hold:

\begin{enumerate}

\item $(M,\cdot ,\tilde{g},E)$ is a weak $\F$-manifold
with Euler vector field $E$ and identity $e:=\tilde{g}^{*}g(E).$
Moreover, $g^{*}\tilde{g}=E\cdot .$

\item If $(g,\tilde{g})$ is a quasi-homogeneous pencil then $e$
is  $\tilde{\nabla}$-flat and $L_{E}(\cdot )=\cdot$ on $TM.$

\end{enumerate}

\end{thm}

We divide the proof into several steps. We begin with the
following Lemma:

\begin{lem}\label{usoara}
Let $(g,\tilde{g})$ be a weak quasi-homogeneous pencil. The
following facts hold:
\begin{enumerate}

\item The multiplication $\lq\lq~\cdot$" on $T^{*}M$ is associative,
commutative and has unity $g(E).$

\item For every $\alpha ,\beta\in{\mathcal E}^{1}(M)$, $g^{*}(\alpha
,\beta )=(\alpha\cdot\beta )(E).$

\end{enumerate}

\end{lem}

\begin{proof}
We prove the first claim: the commutativity $T(u)\cdot
T(v)=T(v)\cdot T(u)$ is equivalent with $T(u)\circ v=T(v)\circ u$
or with $(g(E)\circ u)\circ v=(g(E)\circ v)\circ u$, which follows
from relation (\ref{rightsym}), the metrics $g$ and $\tilde{g}$
being compatible. Using the commutativity $\lq\lq~\cdot$", the
associativity $(u\cdot v)\cdot w=u\cdot (v\cdot w)$ is equivalent
with $(v\cdot u)\cdot w=(v\cdot w)\cdot u$, which is again a
consequence of relation (\ref{rightsym}) and the definition of
\lq\lq\,$\cdot$ ". Relation $T(v)=g(E)\circ v$, for every
$v:=T^{-1}(u)\in T^{*}M$ becomes $g(E)\cdot u=u$, for every $u\in
T^{*}M.$ This clearly means that the multiplication \lq\lq\,$\cdot$" of
$T^{*}M$ has unity $g(E).$

The second claim is an application of the definitions and of
relation (\ref{interpret1}):
\begin{align*}
g^{*}(T(\beta ),\alpha)-(\alpha\circ\beta )(E)
&=g^{*}\left( g(E)\circ\beta ,\alpha\right)-(\alpha\circ\beta )(E)\\
&=g^{*}\left( g(E),\alpha\circ\beta \right)-(\alpha\circ\beta )(E)\\
&=0.
\end{align*}
\end{proof}

\textit{Proof of the Theorem:} The multiplication $\lq\lq~\cdot$"
being associative and commutative on $T^{*}M$, the induced
multiplication $\lq\lq~\cdot$" on $TM$ has the same properties and
has the identity $e=\tilde{g}^{*}g(E).$ Since $(g,\tilde{g})$ is a
weak quasi-homogeneous pencil, $g$ and $\tilde{g}$ are in
particular almost compatible and relation (\ref{interpret2}) with
$\beta$ replaced by $T^{-1}(\beta )$ implies the
$\lq\lq~\cdot$"-invariance of $\tilde{g}.$ An immediate
consequence of the \lq\lq\,$\cdot$"-invariance of $\tilde{g}$ is the
relation $g^{*}\tilde{g}=E\cdot$: indeed, to show this we notice
from the second part of Lemma \ref{usoara} that $g^{*}(\alpha
,\beta )=\tilde{g}^{*}(\alpha\cdot\beta
,E^{\flat})=\tilde{g}^{*}(\alpha\cdot E^{\flat},\beta )$, which
implies that $g^{*}\alpha =\tilde{g}^{*}(\alpha\cdot E^{\flat})$,
or, for $\alpha :=\tilde{g}(X)$, $g^{*}\tilde{g}(X)=E\cdot X.$ In
order to prove that $E$ satisfies the conditions of an Euler
vector field, we notice first that
$L_{E}(\tilde{\nabla})=L_{E}(\nabla )=0$ (since $L_{E}(g)=(1-d)g$,
$L_{E}(\tilde{g})=D\tilde{g}$ and $d\,,D$ are constant), which imply
in particular that $L_{E}(T)=0$. This, together with
$L_{E}(g^{*})=(d-1)g^{*}$, imply that $L_{E}(\cdot )=(d-1)\cdot$
on $T^{*}M$ or, using $L_{E}(\tilde{g})=D\tilde{g}$, that
$L_{E}(\cdot )=(d+D-1)\cdot$ on $TM.$ We can now apply Proposition
\ref{Fman} to prove the weak $\F$-manifold condition (\ref{s}):
the very definition of $\lq\lq~\cdot$", $\lq\lq~\circ$" and the
relation (\ref{cov}) imply that relation (\ref{ec}) is satisfied.
Thus $(M,\cdot ,\tilde{g},E)$ is a weak $\F$-manifold and the
first claim of the Theorem is proved.

The second claim of the Theorem is trivial: if $D=2-d$, then
$L_{E}(\cdot )=\cdot$ on $TM$ and moreover, the very definition of
quasi-homogeneous pencils implies that $e$ is
$\tilde{\nabla}$-flat.

\section{Curvature properties of weak $\F$-manifolds}

In this Section we will prove the first two
vertical 1-1 correspondences in Table 1 of the
introduction. In particular, we show how Dubrovin's correspondence
between flat quasi-homogeneous pencils and Frobenius manifolds
fits into our theory.

We begin with the following simple Lemma on conformal Killing
vector fields on pseudo-Riemannian manifolds.

\begin{lem}\label{confkilling}
Let $(M,\tilde{g})$ be a pseudo-Riemannian manifold and
$E\in{\mathcal X}(M)$ which satisfies
$L_{E}(\tilde{g})=D\tilde{g}$, with $D$ constant. Then
$$
\tilde{g}(\tilde{R}_{Z,X}(E),Y)=g\left(\tilde{\nabla}_{X}(\tilde{\nabla}
E)_{Y},Z\right)
$$
for every $X,Y,Z\in{\mathcal X}(M)$.
\end{lem}

\begin{proof}
We take the covariant derivative with respect to $Z$ of the
equality
$$
\tilde{g}(\tilde{\nabla}_{X}E,Y)+\tilde{g}(\tilde{\nabla}_{Y}E,X)=D\tilde{g}(X,Y),
$$
we use $\tilde{\nabla}(\tilde{g})=0$ and then we take cyclic
permutations of $X$, $Y$, $Z$. We obtain three relations.
Substracting the second and the third relation from the first one
and using the symmetries of pseudo-Riemannian curvature tensors,
we easily obtain the result.
\end{proof}

\noindent Thus if ${\tilde R}=0$, then $E$ can be at most linear
in the flat coordinates.

\begin{thm}\label{flat}
Let $(M,\cdot ,\tilde{g},E)$ be a weak $\F$-manifold as in
Definition \ref{wF}. Then for every $\alpha\in {\mathcal
E}^{1}(M)$ and $X,Y\in{\mathcal X}(M)$, the following relation
holds:
\begin{align*}
R_{E\cdot X,E\cdot Y}(\alpha )-\tilde{R}_{E\cdot X,E\cdot
Y}(\alpha )&=\tilde{\nabla}_{E\cdot X}(\cdot )(T(\alpha )\cdot
(E^{\flat})^{-1}, E^{\flat}\cdot Y^{\flat})
-\tilde{\nabla}_{E\cdot Y}(\cdot )( T(\alpha )\cdot
(E^{\flat})^{-1},
E^{\flat}\cdot X^{\flat})\\
&+\tilde{\nabla}_{E\cdot X}(T)(\alpha )\cdot
Y^{\flat}-\tilde{\nabla}_{E\cdot Y}(T)(\alpha )\cdot X^{\flat}.
\end{align*}
In particular, $M$ is an $\F$-manifold if and only if
$$
R_{E\cdot X,E\cdot Y}(\alpha )-\tilde{R}_{E\cdot X,E\cdot
Y}(\alpha ) =\tilde{\nabla}_{E\cdot X}(T)(\alpha )\cdot
Y^{\flat}-\tilde{\nabla}_{E\cdot Y}(T)(\alpha )\cdot X^{\flat},
$$
for every $\alpha\in {\mathcal E}^{1}(M)$ and $X,Y\in{\mathcal
X}(M)$.
\end{thm}

\begin{proof}
Let $Q (\alpha ):=T(\alpha )\cdot E^{\flat -1}.$ Using relation
(\ref{ec}), we easily see that
\begin{align*}
R_{X,Y}(\alpha )-\tilde{R}_{X,Y}(\alpha )&= [\tilde{\nabla}_{X}(Q
(\alpha ))-Q(\tilde{\nabla}_{X}\alpha
)-Q( Q(\alpha )\cdot X^{\flat}) ]\cdot Y^{\flat}\\
&-[\tilde{\nabla}_{Y}(Q( \alpha ))-Q(\tilde{\nabla}_{Y}\alpha
)-Q( Q(\alpha )\cdot Y^{\flat}) ]\cdot X^{\flat}\\
&+\tilde{\nabla}_{X}(\cdot )(Q(\alpha
),Y^{\flat})-\tilde{\nabla}_{Y}(\cdot )(Q(\alpha ),X^{\flat})\\
\end{align*}
which becomes, after replacing $\alpha$ with $T^{-1}(\alpha )$,
$X$, $Y$ with $E\cdot X$, $E\cdot Y$ respectively, the following
relation:
\begin{align*}
R_{E\cdot X,E\cdot Y}(T^{-1}\alpha )-\tilde{R}_{E\cdot X,E\cdot
Y}(T^{-1}\alpha )&=\tilde{\nabla}_{E\cdot X}(\cdot )(\alpha\cdot
(E^{\flat})^{-1},E^{\flat}\cdot Y^{\flat}) \\
&-\tilde{\nabla}_{E\cdot Y}(\cdot )(\alpha\cdot
(E^{\flat})^{-1},E^{\flat}\cdot X^{\flat})]\\
&+[\tilde{\nabla}_{E\cdot X}(\alpha\cdot (E^{\flat})^{-1})\cdot
E^{\flat}+\tilde{\nabla}_{\tilde{g}^{*}(\alpha )\cdot
X}(E^{\flat})]\cdot Y^{\flat}\\
&-[\tilde{\nabla}_{E\cdot Y}(\alpha\cdot (E^{\flat})^{-1})\cdot
E^{\flat}+\tilde{\nabla}_{\tilde{g}^{*}(\alpha )\cdot
Y}(E^{\flat})]\cdot X^{\flat}\\
&-Q(\tilde{\nabla}_{E\cdot X}(T^{-1}\alpha ))\cdot E^{\flat}\cdot
Y^{\flat}+Q(\tilde{\nabla}_{E\cdot Y}(T^{-1}\alpha ))\cdot
E^{\flat}\cdot X^{\flat}.
\end{align*}
Define now ${\mathcal A}(\alpha ,X):=\tilde{\nabla}_{E\cdot
X}(\alpha\cdot (E^{\flat})^{-1})\cdot E^{\flat}+
\tilde{\nabla}_{\tilde{g}^{*}(\alpha ) \cdot X}(E^{\flat}).$ Using
the fact that $M$ is a weak $\F$-manifold, we easily get
\begin{align*}
{\mathcal A}(\alpha ,X)&=\tilde{\nabla}_{E\cdot X}(\alpha
)-\tilde{\nabla}_{E\cdot X}(\cdot )(E^{\flat},\alpha\cdot
(E^{\flat})^{-1})-\alpha\cdot
(E^{\flat})^{-1}\cdot\tilde{\nabla}_{E\cdot
X}(E^{\flat})+\tilde{\nabla}_{\tilde{g}^{*}(\alpha )\cdot X}(E^{\flat})\\
&=\tilde{\nabla}_{E\cdot X}(\alpha )-\tilde{\nabla}_{E}(\cdot
)(E^{\flat}\cdot X^{\flat},\alpha\cdot
(E^{\flat})^{-1})-\alpha\cdot
(E^{\flat})^{-1}\cdot\tilde{\nabla}_{E\cdot X}(E^{\flat})+\tilde{\nabla}_{\tilde{g}^{*}(\alpha )\cdot X}(E^{\flat})\\
&=\tilde{\nabla}_{E\cdot X}(\alpha
)-\tilde{\nabla}_{E}(\alpha\cdot
X^{\flat})+\tilde{\nabla}_{E}(E^{\flat}\cdot
X^{\flat})\cdot\alpha\cdot (E^{\flat})^{-1}+E^{\flat}\cdot
X^{\flat}\cdot \tilde{\nabla}_{E}(\alpha\cdot (E^{\flat})^{-1})\\
&-\alpha\cdot (E^{\flat})^{-1}\cdot \tilde{\nabla}_{E\cdot
X}(E^{\flat})+\tilde{\nabla}_{\tilde{g}^{*}(\alpha )\cdot
X}(E^{\flat})\\
&=\tilde{\nabla}_{E\cdot X}(\alpha )-\tilde{g}\left(
L_{E}(\tilde{g}^{*}(\alpha )\cdot X)\right) +\tilde{g}\left(
L_{E}(E\cdot X)\right)\cdot \alpha\cdot
(E^{\flat})^{-1}+E^{\flat}\cdot
X^{\flat}\cdot\tilde{\nabla}_{E}(\alpha\cdot
(E^{\flat})^{-1})\\
&=\tilde{\nabla}_{E\cdot X}(\alpha )+X^{\flat}\cdot h(\alpha ),\\
\end{align*}
where in the last equality we have used $L_{E}(\cdot )=k\cdot$ on
$TM$ and $h :=h(\alpha )$ is an expression which depends only on
$\alpha .$ We thus obtain
\begin{align*}
R_{E\cdot X,E\cdot Y}(T^{-1}\alpha )-\tilde{R}_{E\cdot X,E\cdot
Y}(T^{-1}\alpha )=&\tilde{\nabla}_{E\cdot X}(\cdot )(\alpha\cdot
(E^{\flat})^{-1},E^{\flat}\cdot Y^{\flat}) -
\tilde{\nabla}_{E\cdot Y}(\cdot )(\alpha\cdot (E^{\flat})^{-1},E^{\flat}\cdot X^{\flat})\\
&+[\tilde{\nabla}_{E\cdot X}(\alpha )-Q(\tilde{\nabla}_{E\cdot X}(T^{-1}\alpha ))\cdot E^{\flat}]\cdot Y^{\flat}\\
&-[\tilde{\nabla}_{E\cdot Y}(\alpha )-Q(\tilde{\nabla}_{E\cdot Y}(T^{-1}\alpha ))\cdot E^{\flat}]\cdot X^{\flat}\\
\end{align*}
which easily implies the Theorem.
\end{proof}

\begin{cor}(see \cite{dubrovin})
Let $(g,\tilde{g})$ be a quasi-homogeneous pencil on $M$. Suppose
that $\tilde{R}=0.$ Then the corresponding weak $\F$-manifold is a
Frobenius manifold if and only if $R=0.$
\end{cor}

\begin{proof}
This is just a consequence of Theorem \ref{flat} and Lemma
\ref{confkilling}: since $\tilde{R}=0$, $\tilde{\nabla}(T)=0$ and
$$
R_{E\cdot X,E\cdot Y}(\alpha )=\tilde{\nabla}_{E\cdot X}(\cdot
)(T(\alpha )\cdot (E^{\flat})^{-1}, E^{\flat}\cdot Y^{\flat})
-\tilde{\nabla}_{E\cdot Y}(\cdot )( T(\alpha )\cdot
(E^{\flat})^{-1}, E^{\flat}\cdot X^{\flat})
$$
for every $X,Y\in{\mathcal X}(M)$ and $\alpha\in{\mathcal
E}^{1}(M).$ It follows that $R=0$ is equivalent to the
total symmetry of ${\tilde{\nabla}}(\cdot)\,.$
The other conditions of a Frobenius
manifold are clearly satisfied.
\end{proof}

\section{Compatible metrics and submanifolds}

Suppose now that the metrics $g$ and $\tilde{g}$ are compatible.
Let $h$, $\tilde{h}$ be the metrics induced by $g$, $\tilde{g}$ on
a submanifold $N$ of $M$. We assume that $h$, $\tilde{h}$,
$h_{\lambda}$ are non-degenerate (although a theory of
bi-Hamiltonian structures with degenerate metrics may be
formulated \cite{G,ian3}). Let
$A:=\tilde{g}^{*}g$ and $B:=\tilde{h}^{*}h.$\\

\paragraph{\textbf{Notations:}}
We shall use the following conventions:
\begin{enumerate}

\item $TN^{\perp g}$ ($TN^{\perp\tilde{g}}$ respectively)
will denote the orthogonal complement of $TN$ in $TM$, with
respect to the metric $g$ (respectively $\tilde{g}$).

\item For a sub-bundle $V$ of $TM$, $V^{0}$ will denote the
annihilator of $V$ in $T^{*}M.$

\item With respect to the orthogonal decomposition
$TM=TN\oplus TN^{\perp\tilde{g}}$ any tangent vector $X\in TM$
will be written as $X=X^{t}+X^{n}.$
\end{enumerate}

The following Lemma can be easily proved and justifies Definition
\ref{dist}.

\begin{lem}
For every $X\in TN$, $B(X)=A(X)^{t}$. In particular, the following
statements are equivalent:

\begin{enumerate}

\item For every $X\in TN$, $A(X)=B(X).$

\item $A(TN)\subset TN.$

\item The orthogonal complement of $TN$ with respect
to the metrics $g$ and $\tilde{g}$ coincide.

\end{enumerate}

\end{lem}

\begin{defn}\label{dist}
The submanifold $N$ of $M$ is called distinguished if
$A(TN)\subset TN$.
\end{defn}

For now on we shall restrict to the case when the submanifold $N$
is distinguished and we shall denote by $TN^{\perp}$ the
orthogonal complement of $TN$ in $TM$, with respect to the metric
$g$ or $\tilde{g}.$ For $\alpha\in{\mathcal E}^{1}(N)$, we will
denote by $\bar{\alpha}\in{\mathcal E}^{1}(M)$ its extension to
$TM$, which is zero on $TN^{\perp}.$\\

It can be easily verified that the restrictions to $N$ of the
metrics $g_{\lambda}$ generated by $(g,\tilde{g})$ coincide with
the metrics $h_{\lambda}$ generated by $(h,\tilde{h})$. Since
$A=B$ on $TN$, the metrics $h$ and $\tilde{h}$ are almost
compatible. A natural problem which arises is to determine when
they are compatible. For this, let $D$ (respectively $\tilde{D}$)
be the Levi-Civita connections of $h$ and $\tilde{h}$. For every
${\alpha}\in{\mathcal E}^{1}(N)$ and $X\in{\mathcal X}(N)$ we
decompose $\nabla_{X}\bar{\alpha}$ as
\begin{equation}\label{secondf}
\nabla_{X}\bar{\alpha} =\nabla_{X}(\bar{\alpha}
)^{t}+S_{X}(\bar{\alpha})
\end{equation}
according to the decomposition
$$
T^{*}M=(TN^{\perp})^{0}+(TN)^{0}
$$
of $T^{*}M.$ Note that, via the identifications $T^{*}N\cong
(TN^{\perp})^{0}$ and $(TN^{\perp})^{*}\cong (TN)^{0}$,
$D_{X}({\alpha})=\nabla_{X}(\bar{\alpha} )^{t}$ and that the map
$S:TN\times T^{*}N\rightarrow (TN)^{0}$, defined by
$S_{X}({\alpha})=S_{X}(\bar{\alpha})$, is the second fundamental
form of the submanifold $N$ of $(M,g).$ Similar facts hold for
$\tilde{D}$ and the second fundamental form $\tilde{S}$ of the
submanifold $N$ of $(M,\tilde{g}).$

\begin{prop}
The metrics $h$ and $\tilde{h}$ are compatible if and only if for
every $\alpha ,\beta\in (TN^{\perp})^{0}$ and $X,Y\in TN$, the
relation
\begin{equation}\label{submanifold}
\tilde{g}^{*}(\tilde{S}_{X}\alpha -S_{X}\alpha ,
\tilde{S}_{Y}\beta -S_{Y}\beta )=
\tilde{g}^{*}(\tilde{S}_{Y}\alpha -S_{Y}\alpha ,
\tilde{S}_{X}\beta -S_{X}\beta ),
\end{equation}
holds.
\end{prop}

\begin{proof}
This is just a consequence of relation (\ref{conditiecomp}) and of
the decomposition (\ref{secondf}).
\end{proof}

\paragraph{\textbf{Remark:}}
The compatibility condition on the metrics $h$ and $\tilde{h}$ can
also be written in terms of the multiplication $\circ$ from
Definition \ref{multipl}. Denote by
$$
\circ^{t}:T^{*}M\times T^{*}M\rightarrow (TN^{\perp})^{0}
$$
and
$$
\circ^{n}:T^{*}M\times T^{*}M\rightarrow (TN)^{0}
$$
the maps induced by $\circ$ and the orthogonal projections
$T^{*}M\rightarrow (TN^{\perp})^{0}$, $T^{*}M\rightarrow
(TN)^{0}$. Let $\circ_{N}$ be the multiplication associated to the
pair $(h,\tilde{h})$, as in Definition \ref{multipl}. Then, for
every $\alpha ,\beta\in{\mathcal E}^{1}(N)$,
$$
\alpha\circ_{N}\beta =D_{h^{*}\alpha}(\beta
)-\tilde{D}_{h^{*}\alpha}(\beta )
=\nabla_{g^{*}\bar{\alpha}}(\bar{\beta})^{t}-\nabla_{g^{*}\bar{\alpha}}(\bar{\beta})^{t}\\
=\bar{\alpha}\circ^{t}\bar{\beta},
$$
which implies $\alpha\circ_{N}\beta
=\bar{\alpha}\circ^{t}\bar{\beta}.$ It follows that $h$ and
$\tilde{h}$ are compatible if and only if for every
$\bar{\alpha},\bar{\beta},\bar{\gamma}\in (TN^{\perp})^{0}$,
$$
(\bar{\alpha}\circ^{t}\bar{\beta})\circ^{t}\bar{\gamma}=
(\bar{\alpha}\circ^{t}\bar{\gamma})\circ^{t}\bar{\beta}.
$$

\section{Submanifolds of weak $\F$-manifolds}

In this Section we consider a weak $\F$-manifold $(M,\cdot
,\tilde{g},E)$ with $L_{E}(\tilde{g})=D\tilde{g}$ and $L_{E}(\cdot
)=k\cdot .$ Let $N$ be a submanifold of $M$ which satisfies the
following two properties:

\begin{enumerate}

\item For every $X$, $Y$ in $TN$, $X\cdot Y$ belongs to $TN.$

\item The Euler vector field $E$ satisfies
\begin{equation}\label{euler}
E\cdot TN\subset TN.
\end{equation}

\end{enumerate}

We note that any natural submanifold \cite{ian} of an
$\F$-manifold satisfies these conditions, but in our case the
vector field $E$ is not necessarily tangent to $N$ along $N.$ We
shall denote by $TN^{\perp}$ the orthogonal complement of $TN$ in
$TM$, with respect to the metric $\tilde{g}$ and by
$P:TM\rightarrow TN^{\perp}$ the orthogonal projection.

\begin{lem}\label{verysimple}
For every $X\in TN$ and $Y\in TM$, $X\cdot P(Y)=P(X\cdot Y).$
\end{lem}

\begin{proof}
Since $\tilde{g}$ is $\lq\lq~\cdot$"-invariant, the operator
$X\cdot$ of $TM$ is self-adjoint (with respect to the metric
$\tilde{g}$). Since it preserves $TN$, it will preserve
$TN^{\perp}$ as well. The conclusion follows.
\end{proof}

Recall now that on a weak $\F$-manifold we have considered a
second metric $g$, defined by $g^{*}\tilde{g}=E\cdot .$

\begin{prop}
The metrics $g$ and $\tilde{g}$ induce compatible metrics  on $N$.
\end{prop}

\begin{proof}
The almost compatibility is obvious from relation (\ref{euler}).
Let $P^{*}:T^{*}M\rightarrow (TN)^{0}$ be the orthogonal
projection with respect to the metric $\tilde{g}^{*}.$ Relation
(\ref{ec}) together with Lemma \ref{verysimple} imply that, for
every $\alpha ,\beta\in (TN^{\perp})^{0}$, $S_{X}\alpha
-\tilde{S}_{X}\alpha =\frac{1}{2}P^{*}Q(\alpha )\cdot X^{\flat}$,
where, we recall, $Q(\alpha )=\frac{1}{2}\left( (D+k)\alpha
-2\tilde{\nabla}_{\tilde{g}^{*}\alpha}(E^{\flat})\right)\cdot
(E^{\flat})^{-1}.$ It is now obvious that relation
(\ref{submanifold}) is satisfied.

\end{proof}

\section*{Appendix: The semi-simple case}

Recall that we have associated to the pair $(g,\tilde{g})$
an endomorphism $N_{A}$ of $TM$ (see Theorem
\ref{ac}).

\begin{defn}
The pair $(g,\tilde{g})$ is semi-simple if the eigenvalues of the
tensor $N_{A}$ are distinct in every point.
\end{defn}

\begin{thm}\label{semisimple}
Any semi-simple pair of almost compatible metrics is compatible.
\end{thm}

\begin{proof}

If $(g,\tilde{g})$ is a semi-simple pair, we can find coordinates
$(x_{1},\cdots ,x_{n})$ on $M$ such that the tensor
$A:=\tilde{g}^{*}g$ is diagonal:
$\tilde{g}^{*}dx_{i}=\lambda_{i}g^{*}dx_{i}$, for
$i=\overline{1,n}$ and moreover, both $g$ and
$\tilde{g}$ are diagonal:
$g^{*}(dx_{i},dx_{j})=\delta_{ij}{g}^{ii}; \quad
\tilde{g}^{*}(dx_{i},dx_{j})=\delta_{ij}\tilde{g}^{ii}$, for
$i,j=\overline{1,n}$, for some smooth functions $\lambda_{i}\,,g^{ii}$ and $\tilde{g}^{ii}$ .
The almost compatibility
condition implies that the functions $\lambda_{i}$ depend only on
the coordinate $x_{i}$ (see \cite{mok}). Using the formula for the
Cristoffel symbols in a chart and the fact that $\lambda_{i}$
depend only on the coordinate $x_{i}$, we easily see that
$\tilde{\Gamma}^{j}_{ik}-\Gamma^{j}_{ik}=0$
for $i\neq k$ and every $j$. Then
$$
\tilde{\nabla}_{\frac{\partial~}{\partial x_{i}}}dx_{j}-
\nabla_{\frac{\partial~}{\partial x_{i}}}dx_{j}
=-\left(\tilde{\Gamma}_{ii}^{j}-\Gamma_{ii}^{j}\right) dx_{i}
$$
which implies
\begin{align*}
\tilde{g}^{*}\left(\tilde{\nabla}_{\frac{\partial~}{\partial
x_{i}}}dx_{j}-\nabla_{\frac{\partial~}{\partial x_{i}}}dx_{j},
\tilde{\nabla}_{\frac{\partial~}{\partial k_{k}}
}dx_{n}-\nabla_{\frac{\partial~}{\partial x_{k}}}dx_{n}\right) &=
(\tilde{\Gamma}_{ii}^{j}-\Gamma_{ii}^{j})
(\tilde{\Gamma}_{kk}^{n}-\Gamma_{kk}^{n})\tilde{g}^{*}(dx_{i},dx_{k})\\
&=\delta_{ik}\tilde{g}^{ii}
(\tilde{\Gamma}_{ii}^{j}-\Gamma_{ii}^{j})
(\tilde{\Gamma}_{kk}^{n}-\Gamma_{kk}^{n})\,.
\end{align*}
This expression is obviously symmetric in $i$ and $k.$
\end{proof}

\paragraph{}

\medskip
\noindent{\bf Acknowledgements}
\medskip
Financial support was provided by the EPSRC (grant GR/R05093).

\vskip 1cm

\noindent{Authors' addresses:}\par

\bigskip

\begin{tabular}{ll}
Liana David: &(Permanent address): Institute of Mathematics of the
Romanian Academy, \\&Calea Grivitei nr 21, Bucharest, Romania;\\&
e-mail: lili@mail.dnttm.ro\\&\\

&(Present address): Department of Mathematics, University of
Glasgow, \\&Glasgow G12 8QW; \\&e-mail: l.david@maths.gla.ac.uk\\
&\\
Ian A.B. Strachan: &Department of Mathematics, University of
Glasgow, \\&Glasgow G12 8QW; \\&e-mail: i.strachan@maths.gla.ac.uk\\

\end{tabular}

\begin{thebibliography}{99}

\bibitem{DN} Dubrovin, B. and Novikov, S.P., {\em Hydrodynamics of weakly deformed soliton lattices. Differential geometry and Hamiltonian theory,}
Russian Math. Surveys {\bf 44} (1989), no. 6, 35--124.

\bibitem{dubrovin} Dubrovin, B., {\em Flat pencils of metrics and Frobenius manifolds}
in {\sl Integrable systems and algebraic geometry} (Kobe/Kyoto,
1997), 47--72, World Sci. Publishing, River Edge, NJ, 1998.

\bibitem{dubrovin2} Dubrovin, B., {\em On almost duality for Frobenius manifolds}, math.DG/0307374.

\bibitem{F} {Ferapontov, E.V., {\sl Nonlocal Hamiltonian operators of hydrodynamic type:
differential geometry and applications} in {\sl Topics in topology
and mathematical physics}, ed. Novikov, S.P.. Amer. Math. Soc.
Transl. Ser. 2 {\bf 170} (1995).}

\bibitem{ferapontov2} Ferapontov, E.V., {\sl Compatible Poisson brackets of hydrodynamic type},
J. Phys. {\bf A34} (2001), 2377--2388.


\bibitem{G}{Grinberg, N.I., {\sl On Poisson brackets of hydrodynamic type with a
degenerate metric}, Russ. Math. Surv. {\bf 40:4} (1985) 231-232.}

\bibitem{hert} {Hertling, C., {\sl Frobenius manifolds and moduli spaces for singularities},
Cambridge Tracts in Mathematics {\bf 151}, Cambridge University
Press (2002).}

\bibitem{HM} {Hertling, C. and Manin, Yu., {\sl Weak Frobenius manifolds},
Int. Math. Res. Notices {\bf 6} (1999) 277-286.}

\bibitem{Magri} Magri, F., {\em A simple model of the integrable Hamiltomian equation,}
J. Math. Phys. {\bf 19} (1978) 1156-1162.

\bibitem{manin} Manin, Yu., {\em F-manifolds with flat structure and Dubrovin's duality,}
math.DG/0402451.

\bibitem{mok} Mokhov, O.I., {\em Compatible flat metrics},
J. Appl. Math. {\bf 2} (2002), no. 7, 337--370.

\bibitem{ian} Strachan, I.A.B., {\em Frobenius submanifolds}, Journal
of Geometry and Physics, 38 (2001), 285-307.

\bibitem{ian2} Strachan, I.A.B., {\em Frobenius manifolds: natural
submanifolds and induced bi-hamiltonian structures}, Differential
Geometry and its Applications, {\bf 20} (2004), 67-99.

\bibitem{ian3} Strachan, I.A.B., {\em Degenerate Frobenius manifolds and
the bi-Hamiltonian structure of rational Lax equations}, J. Math.
Phys. {\bf 40} (1999), 5058--5079.


\bibitem{V} Vaisman, I., {\em Sur quelques formules du calcul de Ricci global}, Comment. Math. Helv. {\bf 41} (1966/67),
73-87.


\end{thebibliography}
\end{document}